\documentclass[12pt,a4paper]{amsart}

\usepackage{amsmath, amsthm, amssymb}
\usepackage{graphicx}
\usepackage[margin = 2.5cm]{geometry}
\usepackage{mathrsfs}

\newcommand{\x}{\mathbf x}
\newcommand{\y}{\mathbf y}

\usepackage[pdf]{pstricks} 
\usepackage{pstricks-add}
\usepackage{pst-pdf}
\usepackage{pstricks}

\begin{document}

\title{Optimization-based control design techniques and tools}
\author{Pierre Apkarian and Dominikus Noll}  
\thanks{Pierre Apkarian is with ONERA,  2 Av. Ed. Belin, 31055, Toulouse, France 
       {\tt\small Pierre.Apkarian@onera.fr}}        
\thanks{Dominikus Noll is with Universit\'e Paul Sabatier, Institut de Math\'ematiques, 118,
route de Narbonne, 31062 Toulouse, France         
{\tt\small dominikus.noll@math.univ-toulouse.fr}}%
\date{}

\begin{abstract}
Structured output feedback controller synthesis is an exciting recent
concept in modern control design, which bridges between theory and practice
in so far as 
it allows for the first time to apply sophisticated mathematical
design paradigms like $H_\infty$- or  $H_2$-control within
control architectures preferred by practitioners. The new approach to structured $H_\infty$-control,
developed by the authors during the past decade,  is rooted in a change of paradigm in the synthesis algorithms.
Structured design
is no longer be based on solving algebraic Riccati equations or matrix inequalities. 
Instead, optimization-based design techniques are required.
In this essay we
indicate why  structured controller synthesis is central in modern control engineering. We
explain why non-smooth optimization techniques are needed to compute structured control laws, and
we point to software tools which enable practitioners to use these new tools in high technology
applications.

\vspace{.1cm}
\end{abstract}

\maketitle

\section*{Introduction}\label{sect-Motivation} 

In the modern high technology field
control engineers usually face a large variety of concurring  design specifications
such as noise or  gain attenuation in prescribed frequency bands, 
damping, decoupling, constraints on settling- or rise-time, and much else. In addition, 
as plant models 
are generally only approximations of the  true system dynamics,  
control laws have to be robust with respect to uncertainty
in  physical parameters or with regard to un-modeled high frequency phenomena. 
Not surprisingly, such a plethora of 
constraints presents a major challenge for controller tuning,
due not only to the ever growing number of such constraints,
but also because of their very different provenience.

The  steady increase in plant complexity is  
exacerbated by the quest  that regulators should be as simple as possible, easy to
 understand and to tune by practitioners, convenient to hardware implement,  
 and generally
 available at low cost. 
These practical constraints highlight the limited use of Riccati- or LMI-based  controllers, and they are the driving
force for the implementation of  
{\em structured} control architectures. On the other hand this means that 
hand-tuning methods have to be replaced by rigorous algorithmic optimization tools.

\section{Structured Controllers}\label{sect-Struct} 
\begin{figure}[h!]
\hspace*{-.4cm}
\includegraphics[scale=0.45]{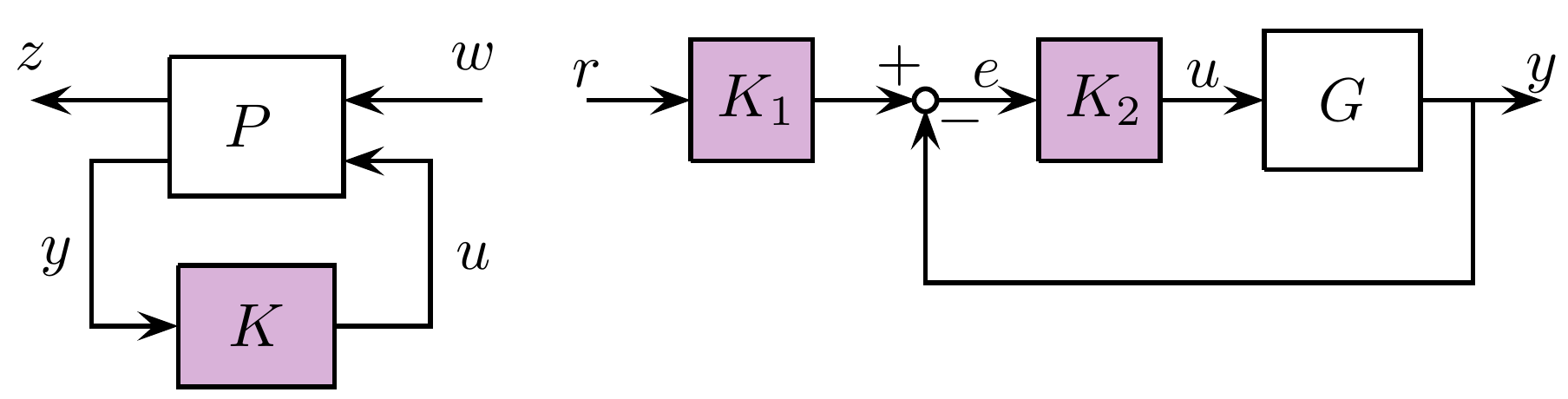}
\caption{Black-box full-order controller $K$ on the left, structured 2-DOF control architecture with
$K={\rm diag}(K_1,K_2)$ on the right.}
\label{fig-BlackBox} 
\end{figure}

Before addressing specific optimization techniques, we recall some basic terminology
for control design problems with structured controllers.  
The  plant model  $P$ is described as 
\begin{eqnarray}
\label{P}
P:\left\{
\begin{array}{llrlrlr}
\dot{x}_P&\!\!=\!\!&Ax_P&\!\!+\!\!&B_1w&\!\!+\!\!&B_2u\\
z&\!\!=\!\!&C_1x_P&\!\!+\!\!&D_{11}w&\!\!+\!\!&D_{12}u\\
y&\!\!=\!\!&C_2x_P&\!\!+\!\!&D_{21}w&\!\!+\!\!&D_{22}u
\end{array}\right.
\end{eqnarray}
where $A$, $B_1$, ... are real matrices of appropriate dimensions,  $x_P\in\mathbb R^{n_P}$ is the state, $u\in\mathbb R^{n_u}$ the control, 
$y\in \mathbb R^{n_y}$ the measured output, $w\in\mathbb R^{n_w}$ the exogenous input,
and $z\in\mathbb R^{n_z}$ the regulated output. Similarly, the sought 
output feedback controller $K$ is described as 
\begin{eqnarray}
\label{K}
K:\left\{
\begin{array}{clrlr}
\dot{x}_K&\!\!\!=\!\!&A_Kx_K&\!\!+\!\!&B_Ky\\
u&\!\!\!=\!\!&C_Kx_K&\!\!+\!\!&D_Ky
\end{array}\right.
\end{eqnarray}
with $x_K\in\mathbb R^{n_K}$, and is called {\em structured} if the (real) matrices $A_K,B_K,C_K,D_K$
depend smoothly on a design parameter $\x\in \mathbb R^n$,
referred to as the vector of tunable parameters. Formally, we have differentiable mappings
\[
A_K = A_K(\x), B_K = B_K(\x), C_K =  C_K(\x), D_K =D_K(\x),
\]
and we abbreviate these by the notation  $K(\x)$ for short to emphasize that the controller is structured
with $\x$ as tunable elements. 
A structured controller synthesis problem
is then an optimization problem of the form
\begin{eqnarray}
\label{program}
\begin{array}{ll}
\mbox{minimize}& \|T_{wz}(P,K(\x))\| \\
\mbox{subject to}& \mbox{$K(\x)$ closed-loop stabilizing}\\
&\mbox{$K(\x)$ \mbox{structured}, \,$\x\in \mathbb R^n$ }
\end{array}
\end{eqnarray} 
where $T_{wz}(P,K)= \mathcal F_\ell(P,K)$ is the lower feedback
connection of (\ref{P}) with (\ref{K}) as in Fig. \ref{fig-BlackBox} (left), also called the Linear Fractional Transformation \cite{ZDG:96}.  The norm $\|\cdot\|$ stands for
the $H_\infty$-norm, the $H_2$-norm, or
any other system norm, while the optimization variable 
$\x\in \mathbb R^n$ regroups the tunable parameters in the design.

Standard examples of structured controllers $K(\x)$ include realizable PIDs, observer-based,
reduced-order, or decentralized controllers, which in state-space are expressed as:
\begin{eqnarray*}
\!\!\!\!\!\!
\begin{array}{c}
\left[ \begin{array}{cc|c}  0&0&1\\0&-1/\tau&-k_D/\tau\\ \hline k_I&1/\tau&k_P +k_D/\tau \end{array} \right], 
 \left[  \begin{array}{c|c}  A-B_2K_c-K_fC_2&K_f\\ \hline -K_c&0 \end{array} \right] ,\\
 \left[  \begin{array}{c|c} A_K&B_K\\ \hline C_K&D_K \end{array} \right] ,\\
 \left[ \begin{array}{c|c} {\rm diag}(A_{K_1},\ldots,A_{K_q})& {\rm diag}(B_{K_1},\ldots,B_{K_q})\\ \hline
 {\rm diag}(C_{K_1},\ldots,C_{K_q})& {\rm diag}(D_{K_1},\ldots,D_{K_q})\\
\end{array}\right].
\end{array}
\end{eqnarray*}
In the case of a  PID the tunable parameters are $\x = (\tau,k_P,k_I,k_D)$,
for observer-based controllers $\x$ regroups the estimator and state-feedback gains
 $(K_f,K_c)$, for reduced order controllers $n_K < n_P$
the tunable parameters
$\x$ are the $n_K^2 + n_Kn_y + n_Kn_u + n_yn_u$ unknown entries in
$(A_K,B_K,C_K,D_K)$, and in the decentralized form
$\x$ regroups the unknown entries in  $A_{K1},\dots,D_{Kq}$.
In contrast,  full-order controllers have the maximum number
$N = n_P^2 + n_Pn_y+n_Pn_u+n_yn_u$ of degrees of freedom and are referred to
as unstructured or as {\em black-box} controllers.
\\

More sophisticated controller  structures  $K(\x)$ arise form architectures
like for instance a 2-DOF control arrangement with feedback block $K_2$ and
a set-point filter $K_1$ as in Fig. \ref{fig-BlackBox} (right).
Suppose 
$K_1$ is  the $1$st-order filter $K_1(s)=a/(s+a)$ and $K_2$ the PI feedback $K_2(s) = k_{P}+ k_{I}/s$.
Then  
the transfer  $T_{ry}$ from $r$ to $ y$ can be represented as the feedback connection
of $P$ and $K(\x,s)$ with 
\begin{eqnarray*}
P:=\left[  \begin{array}{c|ccc} A&0&0&B\\\hline C&0& 0&D\\0&I&0&0\\ -C&0&I&-D \end{array}\right], \\
K(\x,s):=\left[  \begin{array}{cc} K_1(a,s)&0\\0&K_2(s,k_P,k_I)\\
 \end{array} \right]\,,
\end{eqnarray*}
and gathering tunable elements in  
$\x = (a,k_P,k_I)$.

\begin{figure}[h!]
\includegraphics[scale=0.44]{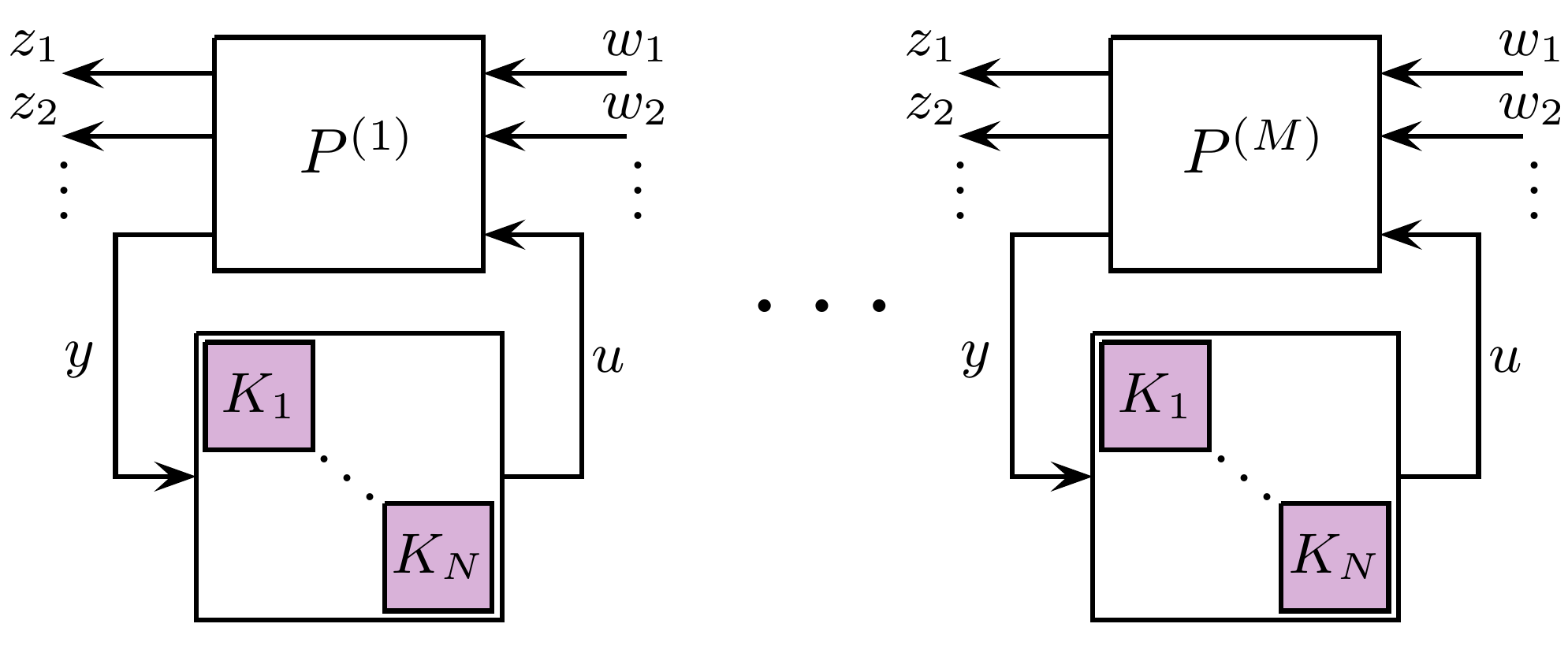}
\caption{Synthesis of $K={\rm diag}(K_1,\dots,K_N)$ against multiple requirements or models 
$P^{(1)},\dots,P^{(M)}$. Each $K_i(\x)$ can  be structured. \label{fig-standardForm2}}
\end{figure}

In much the same way arbitrary multi-loop  interconnections of fixed-model elements 
with tunable controller blocks $K_i(\x)$ can be re-arranged as 
in Fig. \ref{fig-standardForm2}, so that  $K(\x)$ captures all tunable blocks in a decentralized
structure general enough to cover most engineering applications.

The structure concept is equally useful
 to deal with
the second central challenge in control design:  {\em system uncertainty}. 
The latter may be handled
with  $\mu$-synthesis techniques \cite{SD1991}  if a parametric uncertain model
is available. A less ambitious but often more practical alternative
consists in  optimizing  the structured controller $K(\x)$
against a finite set of plants $P^{(1)},\dots,P^{(M)}$ representing model variations due to uncertainty, aging,
sensor and actuator breakdown, un-modeled dynamics, in tandem 
with the robustness and performance specifications. 
This is again formally covered by
Fig. \ref{fig-standardForm2} and
leads to a multi-objective constrained optimization problem
of the form

\begin{eqnarray} \label{eq-objCons}
\begin{array}{ll}
\mbox{minimize} & f(\x)=\max\limits_{k\in {\rm SOFT}, i\in I_k} \| T_{w_i z_i }^{(k)} (K(\x))\| \\
\mbox{subject to} 
&g(\x)= \max\limits_{k\in {\rm HARD}, j\in J_k}  \|T_{w_j z_j}^{(k)} (K(\x)) \|  \leq 1 \\
&K(\x) \mbox{ structured and  stabilizing}  \\
&\x \in \mathbb R^n
\end{array}
\end{eqnarray}
where  $T_{w_iz_i }^{(k)}$ denotes the $i$th
closed-loop robustness or performance channel $w_i \to z_i$ for the $k$-th plant model $P^{(k)}$.
SOFT and HARD denote index sets taken over a finite set of specifications, say in  $\left\{1,\ldots,M\right\}$. The rationale of (\ref{eq-objCons}) is to minimize the
worst-case cost of the soft constraints $\| T_{w_iz_i}^{(k)} \|$, $k\in\,$SOFT,  while enforcing
the hard constraints $\| T_{w_jz_j}^{(k)} \| \leq 1$, $k\in \,$HARD, 
which prevail over soft ones and are mandatory. In addition to local optimization
(\ref{eq-objCons}), the problem can undergo a global optimization step in order to 
prove global stability and performance of the design, see \cite{RNA:2016,AN:2015,ANRsiam:2015}.

\section{Optimization Techniques Over the Years}\label{sect-Opt}
During the late 1990s the necessity to develop
design techniques  for  structured regulators
$K(\x)$ was  recognized \cite{FAN:99}, and
the limitations of synthesis methods based on algebraic Riccati equations 
or linear matrix inequalities (LMIs)  became evident, as
these techniques cannot provide structured controllers needed in practice. 
The lack of appropriate synthesis
techniques for structured $K(\x)$  
led to the unfortunate situation, where  sophisticated approaches like
the $H_\infty$ paradigm developed by academia since the 1980s could not be brought to work 
for the design of those controller structures $K(\x)$
preferred by practitioners. Design engineers had to continue to rely on heuristic and  ad-hoc tuning techniques,
with only limited scope and reliability. 
As an example:  post-processing to reduce a black-box controller to a practical size
is prone to failure.
It may  at best be considered a fill-in for a rigorous design method which directly computes a reduced-order
controller. 
Similarly, hand-tuning of the parameters $\x$ remains a puzzling task because of the loop interactions, and
fails as soon as  complexity  increases. 

In the late 1990s and early 2000s, a change of methods was observed.
Structured $H_2$- and $H_\infty$-synthesis 
problems (\ref{program}) were addressed by 
bilinear matrix inequality (BMI) optimization, which used local optimization
techniques based on the augmented Lagrangian method \cite{FAN:99,NA:00,KS:2003,noll:07}, 
sequential semidefinite programming methods \cite{FNA:00,ANTT:00},
and non-smooth methods for BMIs \cite{NPA:2008,LO:00}. However, these techniques were
based on the bounded real lemma or similar matrix inequalities,
and were therefore of limited success due to the presence of Lyapunov variables, i.e. matrix-valued unknowns, 
whose dimension grows quadratically in $n_P+n_K$ and represents the bottleneck of that approach.

The epoch-making change occurs with the 
introduction of non-smooth optimization techniques \cite{na:03_1,ANtac:05,ANmultiband:06,an_disk:05} to programs
 (\ref{program}) and (\ref{eq-objCons}).   Today non-smooth methods have superseded 
matrix inequality-based techniques and may be considered the state-of-art 
as far as realistic applications are concerned. The transition took almost a decade.

Alternative control-related local optimization techniques and heuristics include the gradient sampling technique 
of \cite{BLO:gradient},  and other derivative-free optimization as in \cite{KLT:2003,an:05}, 
particle swarm optimization, see \cite{ONMFMNAK:2008} and references therein,  
and also evolutionary computation techniques \cite{Lie:2001}. All these  
classes do not exploit derivative information and rely on function evaluations only. 
They are therefore applicable to a broad variety of problems including those where 
function values arise from complex numerical simulations. 
The combinatorial nature of these techniques, however,  limits their use to small problems
with  a few tens of variable.  
More significantly, these methods often lack a solid convergence theory. 
In contrast, as we have demonstrated over recent years, \cite{ANtac:05,NPR:2008,ANR:2015,ANR:17}
specialized non-smooth techniques are highly efficient in practice, are based on a sophisticated convergence theory, 
capable of solving medium size problems in a matter of seconds, and are still operational for
large size problems with several hundreds of states.

\section{Non-smooth optimization techniques}\label{sect-nonsmooth}
The benefit of the non-smooth casts (\ref{program}) and (\ref{eq-objCons})  lies in the possibility to avoid searching for Lyapunov variables, a major advantage
as their number $(n_P+n_K)^2/2$ usually largely dominates $n$, the number of true
decision parameters $\x$.
Lyapunov variables do still occur implicitly in the function evaluation  procedures, but this has 
no harmful effect for systems up to several hundred states. 
In abstract terms, a non-smooth optimization program
has the form
\begin{eqnarray}\label{eq-cast}
\begin{array}{ll}
\mbox{minimize} & f(\x)\\
\mbox{subject to}& g(\x)\leq 0 \\
&\x \in\mathbb R^n
\end{array}
\end{eqnarray}
where $f,g:\mathbb R^n\to \mathbb R$ are locally Lipschitz functions and are easily identified from the cast in (\ref{eq-objCons}). 

In the realm of convex optimization, non-smooth programs are conveniently addressed by so-called
bundle methods, introduced in the late 1970s by Lemar\'echal \cite{Lem:75}.  Bundle
methods are used to solve difficult problems
in integer programming or in stochastic optimization via Lagrangian
relaxation. Extensions of the bundling technique to non-convex problems
like (\ref{program}) or (\ref{eq-objCons}) were first developed in \cite{ANtac:05,ANmultiband:06,an_disk:05,ANP:2007,NPA:2008},
and in  more abstract form,  in \cite{NPR:2008}. Recently, we also extended bundle techniques to the trust-region framework
\cite{ANR:2015,ANR:17,AN:18}, which leads to the first extension of the classical trust-region method to
non-differential optimization supported by a valid convergence theory.

\begin{figure}[h!]
\centering
\includegraphics[width=0.65\textwidth]{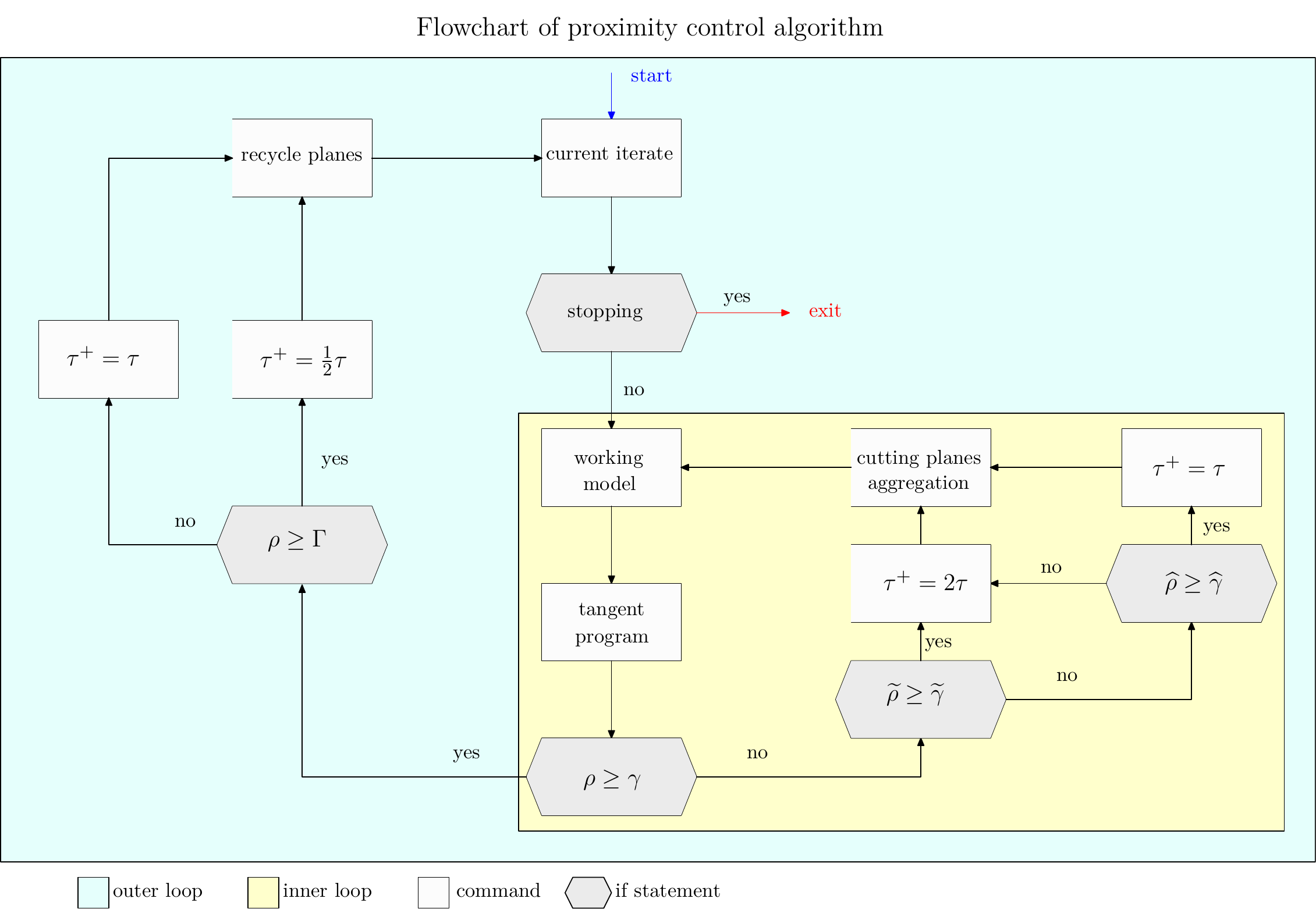}
\caption{Flow chart of proximity control bundle algorithm \label{algo}}
\end{figure}

Fig. \ref{algo}
shows a schematic view of a non-convex bundle method consisting of a descent-step generating
inner loop (yellow block) comparable to a line search in smooth optimization, embedded into the
outer loop (blue box),  where serious iterates are processed, stopping criteria are applied,  
and the acceptance rules of traditional trust region techniques  is assured. At the core of
the interaction between inner and outer loop is the management of the proximity control 
parameter $\tau$, which  governs the stepsize $\|\x-\y^k\|$ between trial steps $\y^k$ at the current
serious iterate $\x$.  Similar to the management of a trust region radius or of the stepsize
in a linesearch,  proximity control allows to force shorter trial steps if agreement of the local  model
with the true objective function is poor, and allows larger steps if agreement is satisfactory.

Oracle-based bundle methods traditionally assure global convergence in the sense
of subsequences under the sole hypothesis that for every trial
point $\x$ the function value $f(\x)$ and one Clarke subgradient $\phi\in \partial f(\x)$
are provided. In automatic control applications it is as a rule possible to
provide more specific information, which may be exploited 
to speed up  convergence \cite{ANtac:05}.

Computing function value and gradients of the
$H_2$-norm $f(\x)=\|T_{wz}\left(P,K(\x)\right)\|_2$
requires essentially the solution of two Lyapunov equations of size $n_P+n_K$,
see \cite{ANR:2007}.      
For the $H_\infty$-norm, 
$f(\x)=\|T_{wz}\left( P,K(\x)  \right)\|_\infty$,  function evaluation is based
on the Hamiltonian algorithm of \cite{BSV:2012,BBK:89}. The Hamiltonian matrix is of size $n_P+n_K$, so that
function evaluations may be costly for very large plant state dimension ($n_P>500$), even though 
the number of outer loop iterations
of the bundle algorithm is not affected by a large $n_P$ and generally relates to $n$, the dimension of $\x$.
The additional cost for subgradient computation for large $n_P$
is relatively cheap as it relies on linear algebra \cite{ANtac:05}. Function and subgradient evaluations for $H_\infty$ and $H_2$
norms are typically obtained in $\mathcal{O} \left((n_P+n_K)^3\right)$ flops. 

\section{Computational Tools}
Our
non-smooth optimization methods  became 
available to the  engineering community
since 2010 via the MATLAB Robust Control Toolbox  \cite{RCT2012b,GA2011a}. 
Routines  {\tt hinfstruct}, {\tt looptune}$\,$ and  {\tt systune}$\,$   are versatile enough 
to define and combine tunable blocks $K_i(\x)$, 
to build and aggregate multiple models and design requirements on $T_{wz}^{(k)}$  of different nature,
and to provide suitable validation tools. Their implementation was carried out in 
cooperation with P. Gahinet (MathWorks). 
These routines   further exploit the structure of  problem (\ref{eq-objCons}) to enhance efficiency, 
see \cite{ANmultiband:06} and \cite{ANtac:05}. 

It should be mentioned that
design problems with multiple hard constraints 
are inherently complex
and generally NP-hard, so that 
exhaustive methods fail even for small to medium size problems. The principled 
decision made in \cite{ANtac:05}, and reflected in the MATLAB tools,  
is  to rely on local optimization techniques instead.
This leads to weaker convergence certificates, but has the advantage to work successfully in practice. 
In the same vein, in (\ref{eq-objCons}) it is preferable to rely on a mixture of soft and hard requirements,
for instance, by the use of  exact penalty functions \cite{na:03_1}.  Key features 
implemented 
in the mentioned MATLAB routines are discussed in \cite{A:2013,GA2011a,ANmultiband:06}.

\section{Applications}
Design of a feedback regulator is an interactive process, in which
tools like {\tt systune}, {\tt looptune} $\,$ or {\tt hinfstruct} $\,$ support the 
designer in various ways. In this section we illustrate their enormous potential  
by 
showing that even infinite-dimensional 
systems may be successfully addressed by these tools. For a plethora of design examples 
for real-rational systems including parametric and complex dynamic uncertainty we refer to   
\cite{RNA:2016,ANR:2015,AN:2015,AN:18,ANR:17}.
For recent applications of our tools in real-world applications see also \cite{falcoz}, where it is in particular explained how
{\tt hinfstruct} helped in 2014 to save the Rosetta mission. Another important application of {\tt hinfstruct} is the design
of the atmospheric flight pilot for the ARIANE VI launcher
 by the Ariane Group 
\cite{GDC2017}.

\subsection{Illustrative example}
We discuss boundary control of a wave equation with anti-stable
damping, 
\begin{align}
\label{wave}
x_{tt}(\xi,t) &= x_{\xi\xi}(\xi,t), \;\; t \geq 0, \xi \in [0,1] \notag \\
x_\xi(0,t) &= -q x_t(0,t), \; \; q > 0, q\not=1 \\
x_\xi(1,t) &= u(t).\notag
\end{align}
where notations $x_t$, $x_{\xi}$ are partial derivatives of $x$ with respect to time and space, respectively. 
In (\ref{wave}), $x(\cdot,t),x_t(\cdot,t)$ is the state, the control 
applied at the boundary $\xi=1$ is $u(t)$, and we assume that the measured
outputs are 
\begin{equation}
\label{outputs}
y_1(t) = x(0,t), y_2(t)=x(1,t),  y_3(t) = x_t(1,t).
\end{equation}
The system has been discussed previously
in \cite{smyshlyaev_krstic:09,Fridman:2015,bresch_krstic:14}  and has been proposed for the control of slip-strick vibrations in drilling devices
\cite{saldivar:11}. Here measurements $y_1,y_2$ correspond to the angular positions of the drill string at the top and bottom level,
$y_3$ measures angular speed at the top level, while control corresponds to a reference velocity at the top. The friction characteristics at the bottom level are
characterized by the parameter $q$, and the control objective is to maintain a constant angular velocity at the bottom.

Similar models have been used to control pressure fields in duct combustion dynamics,
see \cite{DeQueiroz_Rahn:02}.
The challenge in (\ref{wave}), (\ref{outputs})  is to
design implementable controllers despite the use of an infinite-dimensional system model. 

The transfer function of (\ref{wave}) is obtained  from:
\[
G(\xi,s) =\frac{x(\xi,s)}{u(s)} = \frac{1}{s} \cdot \frac{(1-q)e^{s\xi}+(1+q)e^{-s\xi}}{(1-q)e^s - (1+q)e^{-s}},
\]
which in view of (\ref{outputs}) 
leads to $G(s) =[G_1;G_2;G_3]= [G(0,s); G(1,s);sG(1,s)]$. 

\begin{figure}[ht!]
\includegraphics[scale=1.2]{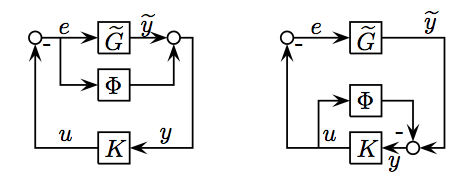}
\caption{Stability of the closed-loop $(\widetilde{G}+\Phi,K)$ is equivalent to stability of the closed-loop $(\widetilde{G},{\tt feedback}(K,\Phi))$.
 \label{swap1}}
\end{figure}

Putting $G$ in feedback with the controller $K_0=[ 0\; 0\; 1]$ leads to
$\widehat{G} = G/(1+G_3)$, where
\begin{equation}
\label{dec}
\widehat{G}(s)= \begin{bmatrix}  \frac{1}{s(1-q)}  \vspace{.2cm}    \\ \frac{1+Q}{2s} 
 \vspace{.2cm}\\ \frac{1}{2}\end{bmatrix}
+ \begin{bmatrix}  -\frac{1-e^{-s}}{s(1-q)} \vspace{.2cm}    \\ -\frac{Q(1-e^{-2s})}{2s} \\ \vspace{.2cm}\frac{Q}{2} e^{-2s}\end{bmatrix} =:\widetilde{G}(s)+\Phi(s).
\end{equation}
Here $\widetilde{G}$ is real-rational and unstable, while $\Phi$ is stable but infinite dimensional.
Now we use the fact that  stability of the closed loop $(\widetilde{G}+\Phi,K)$ is equivalent to
stability of  the loop  $(\widetilde{G}, {\tt feedback}(K,\Phi))$ upon defining 
$
{\tt feedback}(M,N):= M(I+NM)^{-1}\,.
$
The loop transformation is explained  in Fig. \ref{swap1}, see also \cite{Moelja_Meinsma:03}.

Using (\ref{dec})  we construct a finite-dimensional structured
controller $\widetilde{K}=\widetilde{K}(\x)$ which stabilizes $\widetilde{G}$. The controller $K$ stabilizing 
$\widehat{G}$ in (\ref{dec})  is then recovered from $\widetilde{K}$ through the
equation $\widetilde{K} = {\tt feedback}(K,\Phi)$, which when inverted gives
$K = {\tt feedback}(\widetilde{K},-\Phi)$. The overall controller for (\ref{wave}) is 
$K^* = K_0 + K$,  and since along with $K$ only delays appear in $\Phi$, the controller $K^*$ is implementable. 

\begin{figure}[ht!]
\centering
\includegraphics[height=5.2cm]{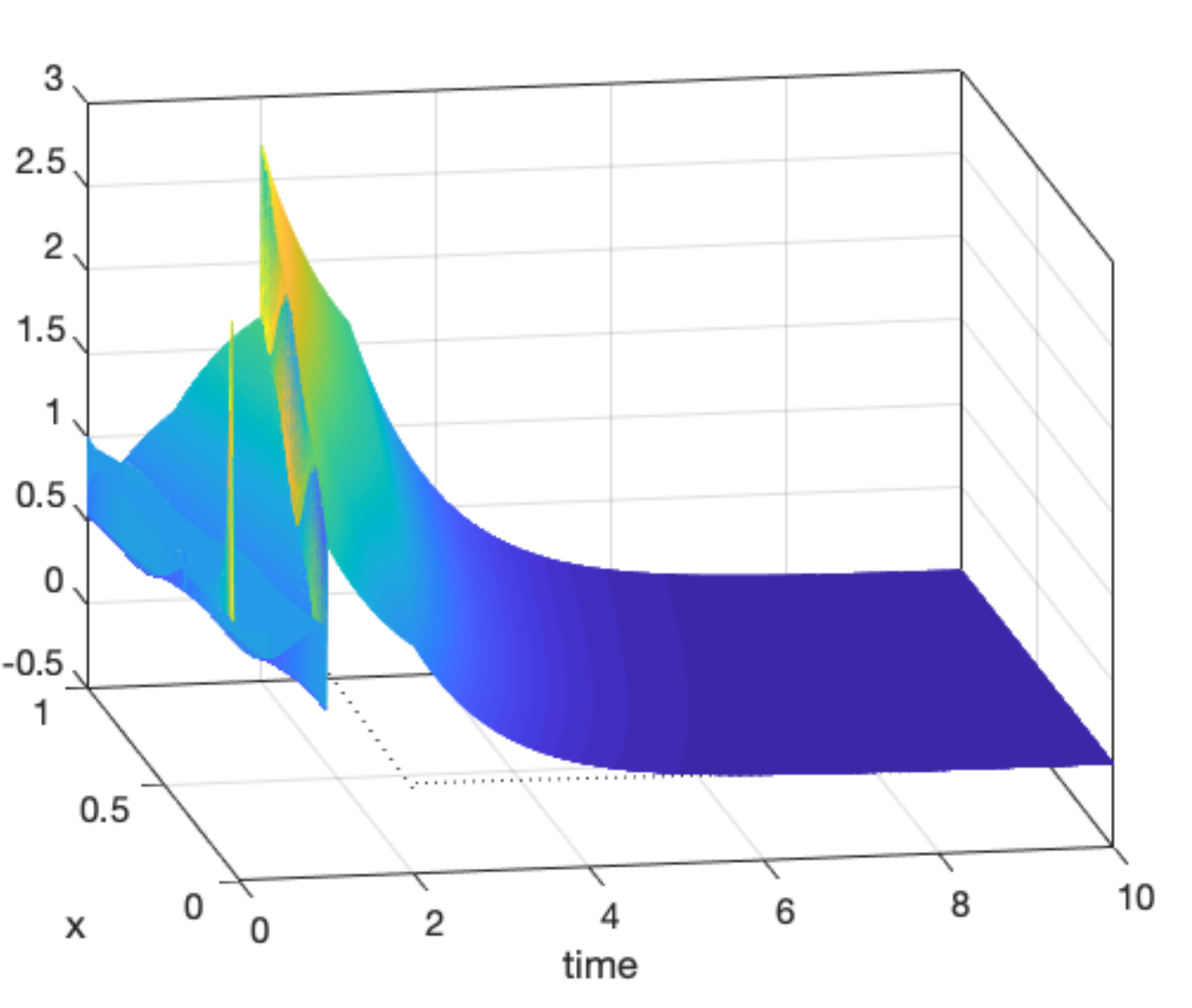}
\includegraphics[height=5.2cm]{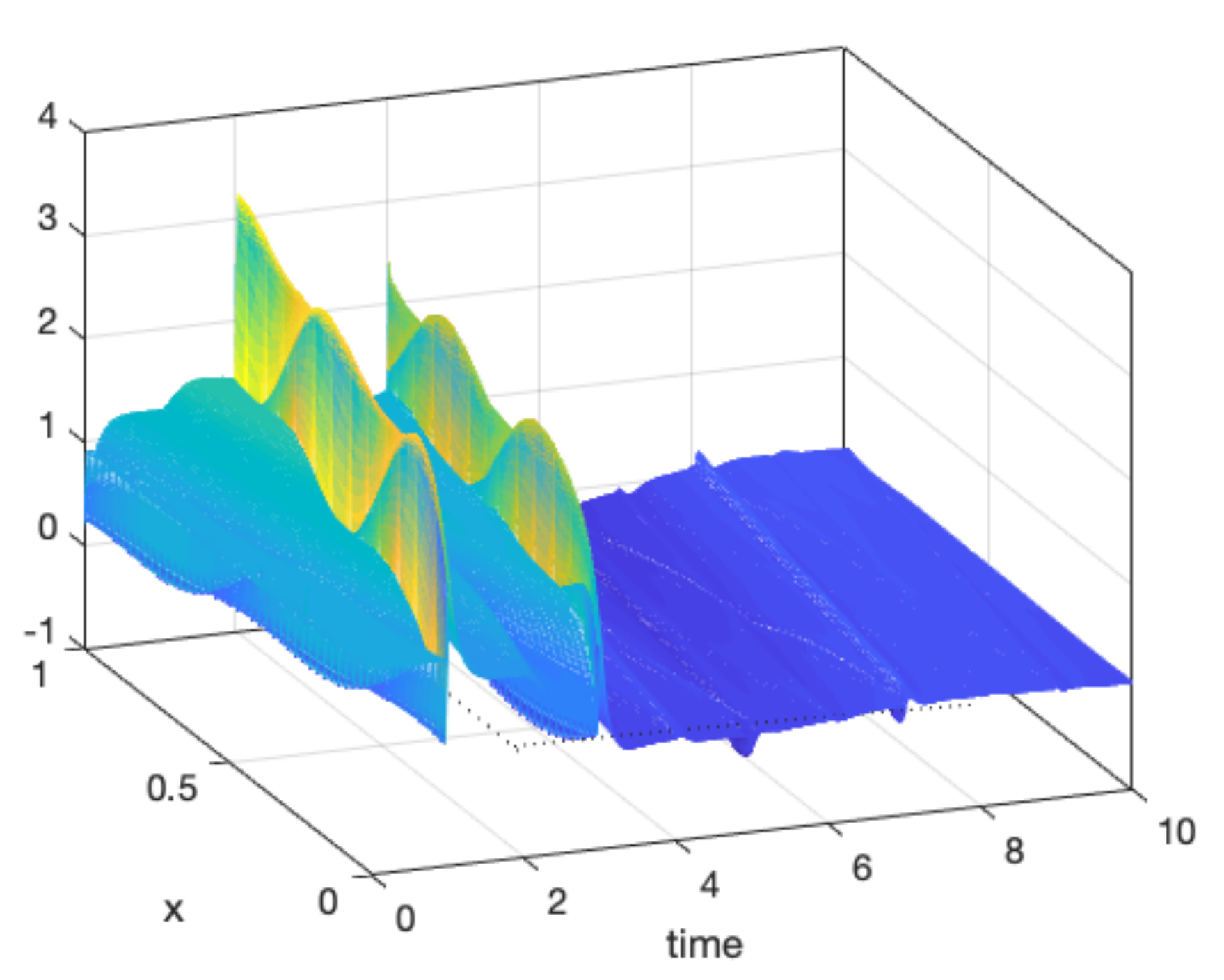}
\caption{Wave equation. Simulations 
for $K$ obtained by backstepping control (top) \cite{bresch_krstic:14}
and  $K^*=K_0+K$ obtained by optimizing ${\tt feedback}(\widetilde{G},\widetilde{K})$ via {\tt systune} $\,$(bottom).
Both controllers are $\infty$-dimensional, but implementable.
\label{Krstic_controller}}
\end{figure}

Construction of $\widetilde{K}$ uses {\tt systune} $\,$with pole placement
via {\tt TuningGoal.Poles},  imposing that closed-loop poles have a minimum decay of 0.9, minimum damping of 0.9, and a maximum frequency of 4.0. The controller structure is chosen as
static, so that $\x \in \mathbb R^3$. A simulation with $K^*$
is shown in Fig. \ref{Krstic_controller} (bottom) and some acceleration over the   backstepping controller from \cite{bresch_krstic:14}
(top)  is observed.

\subsection{Gain-scheduling control}
Our last study is when the parameter $q\geq 0$
is uncertain or allowed to vary in time with sufficiently slow variations as in \cite{shamma90}. We assume that
a nominal  $q_0>0$ and an uncertain interval $[\underline{q},\overline{q}]$
with $q_0 \in (\underline{q},\overline{q})$ and $1\not\in [\underline{q},\overline{q}]$ are given.

The following scheduling
scenarios, all leading to implementable controllers, are possible.
(a) Computing a nominal controller $\widetilde{K}$ at $q_0$ as before, and 
scheduling through $\Phi(q)$, which depends explicitly on $q$, so that
$K^{(1)}(q) = K_0+ {\tt feedback}(\widetilde{K},-\Phi(q))$. (b) Computing 
$\widetilde{K}(q)$ which depends on $q$, 
and using $K^{(2)}(q) = K_0+{\tt feedback}(\widetilde{K}(q),-\Phi(q))$.

While (a) uses (\ref{program}) based on \cite{ANtac:05,an_disk:05} and
available in {\tt systune},  
we show that one can also apply (\ref{program}) to case (b). We use 
Fig. \ref{swap1}  to work in the finite-dimensional system $(\widetilde{G}(q),\widetilde{K}(q))$, where plant and controller
now depend on
$q$, which is a parameter-varying  design.

For that we have to decide on
a parametric form of the controller $\widetilde{K}(q)$, which we choose as
\[
\widetilde{K}(q,\x) = \widetilde{K}(q_0) + (q-q_0) \widetilde{K}_1(\x) + (q-q_0)^2 \widetilde{K}_2(\x),
\]
and 
where we adopted the simple static form
 $\widetilde{K}_1(\x) =[\x_1\; \x_2\; \x_3], \widetilde{K}_2=[\x_4\;\x_5\; \x_6]$, featuring a total of 6 tunable parameters.
The nominal $\widetilde{K}(q_0)$ is obtained via (\ref{program}) as above.
For $q_0 = 3$ this leads to 
$\widetilde{K}(q_0)=[ -1.049 \; -1.049 \; -0.05402]$, computed
via {\tt systune}.

\begin{figure}[ht!]
\includegraphics[height=5.2cm]{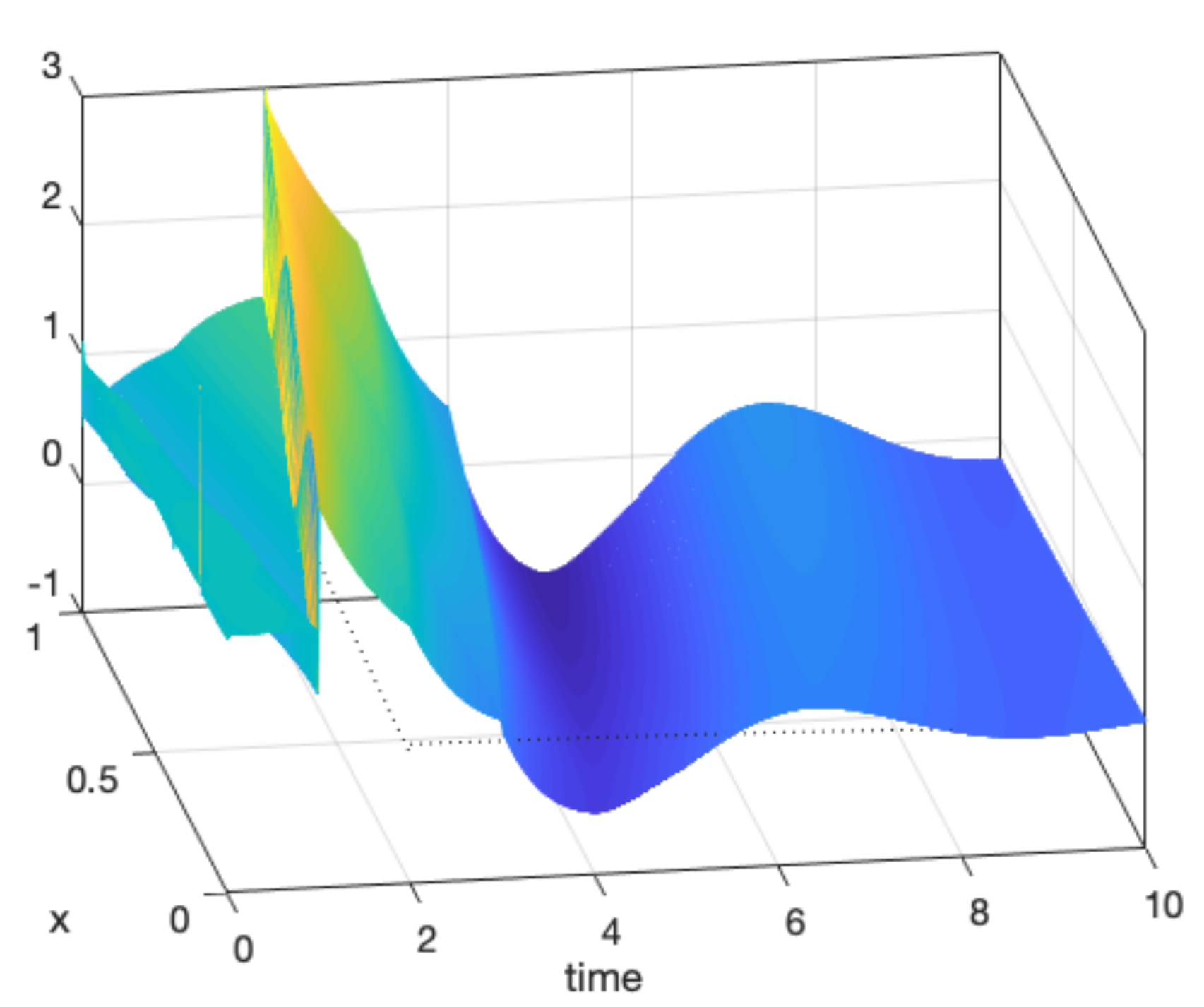}
\includegraphics[height=5.2cm]{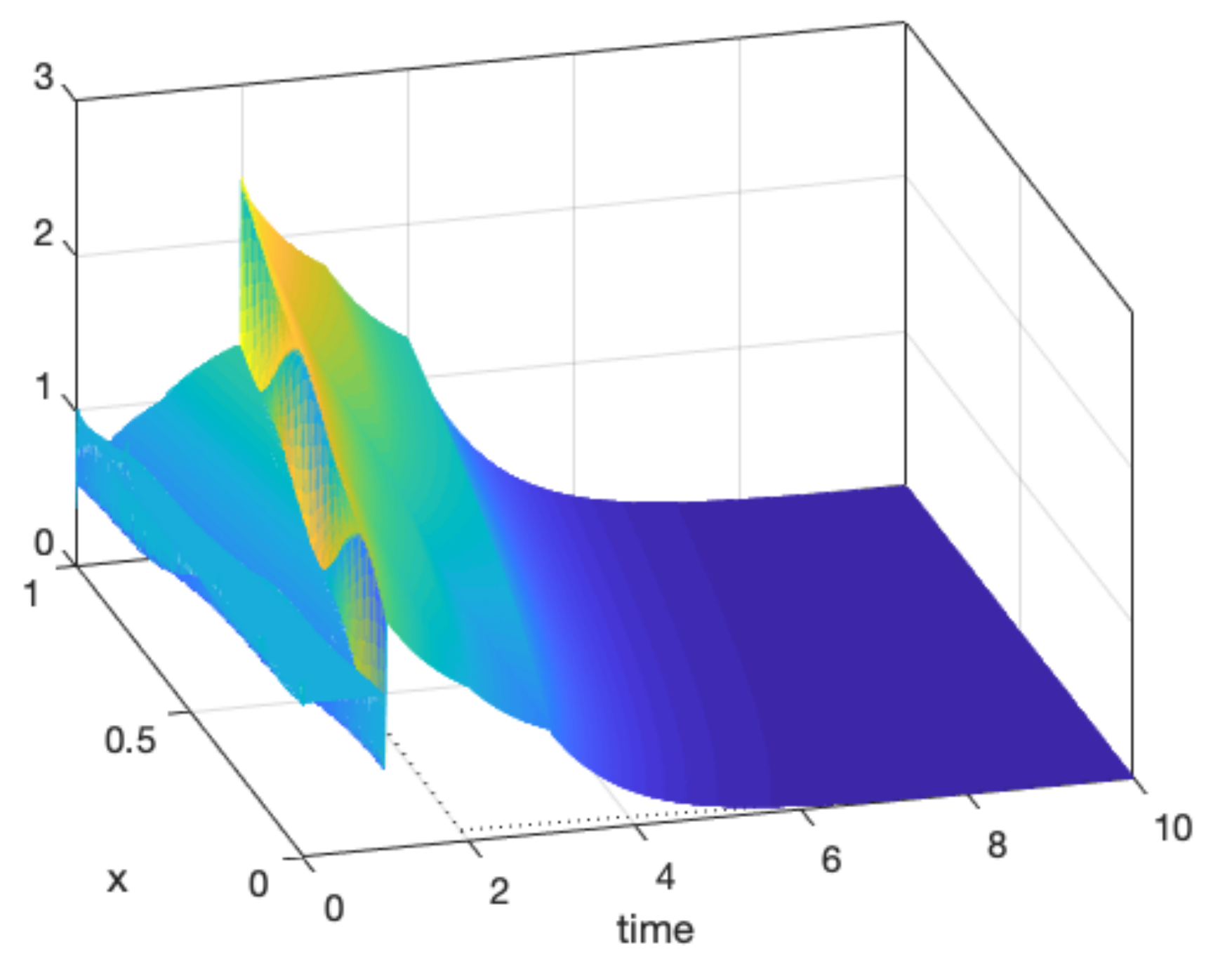}
\includegraphics[height=5.2cm]{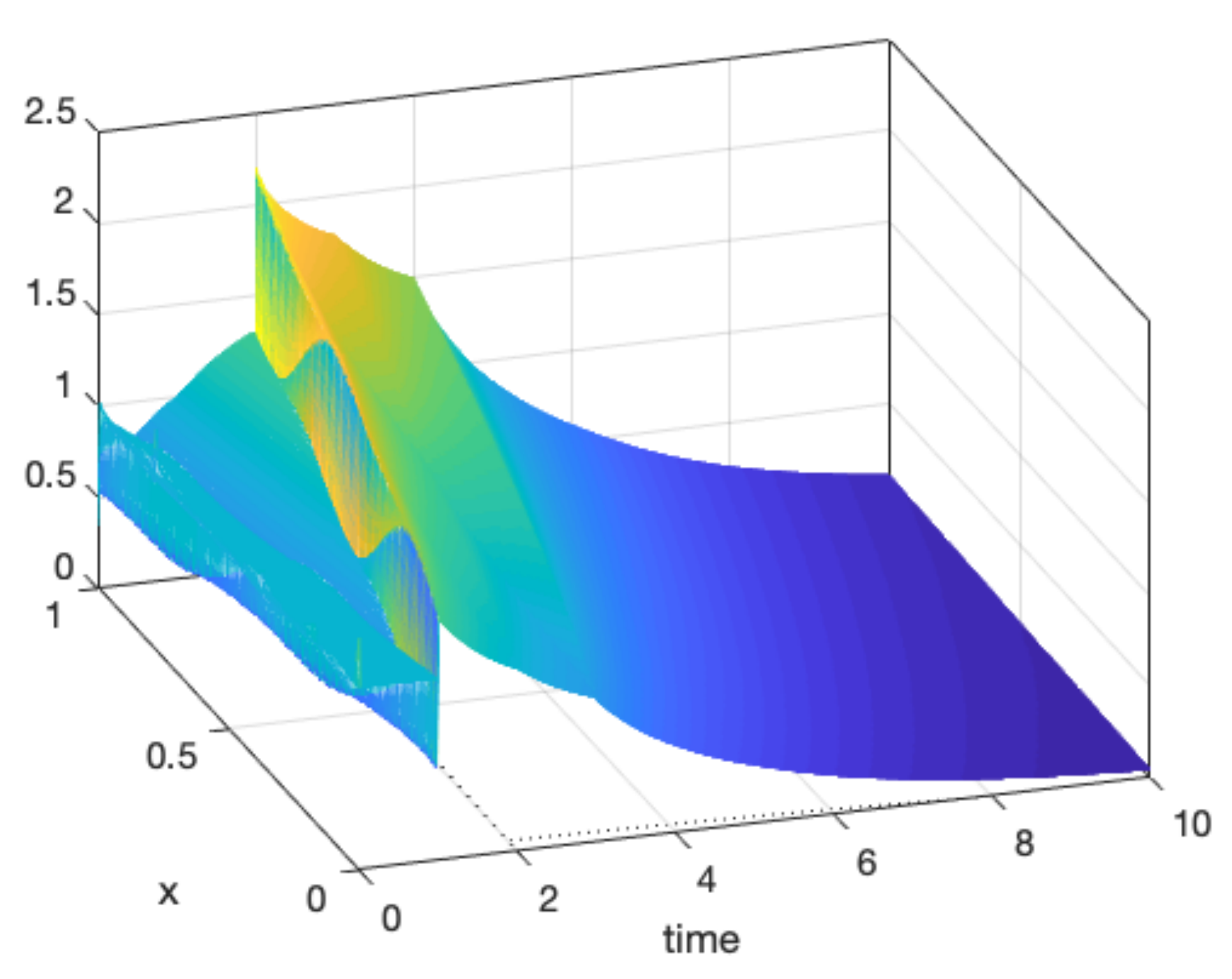}
\caption{Synthesis at nominal $q_0=3$. Simulations of nominal $K=K_0+{\tt feedback}(\widetilde{K},\Phi(3))$  for $q=2,3,4$. Nominal controller is robustly
stable over $[\underline{q},\overline{q}]$.
\label{nominal_robust}}
\end{figure}

\begin{figure}[ht!]
\includegraphics[height=5.2cm]{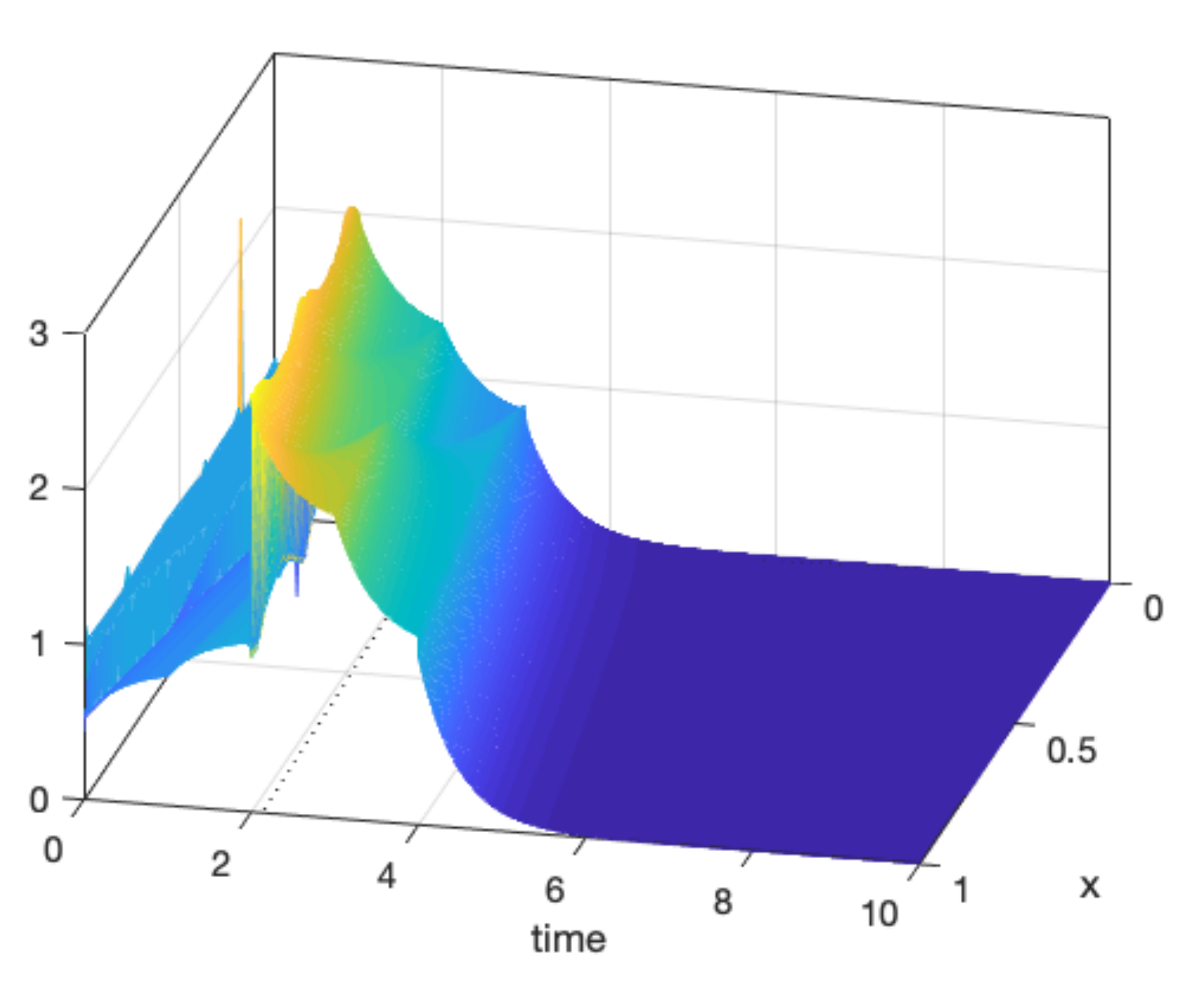}
\includegraphics[height=5.2cm]{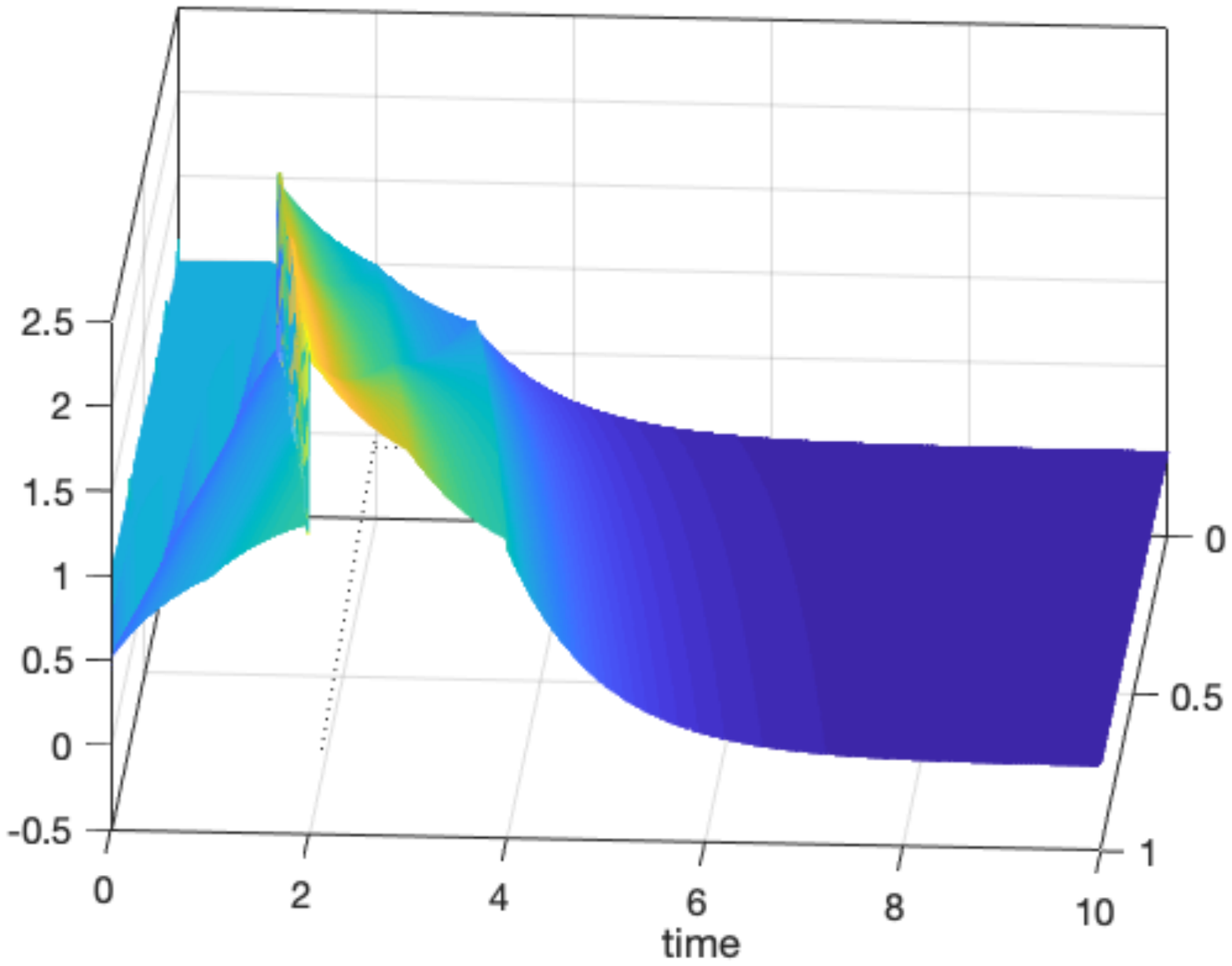}
\includegraphics[height=5.2cm]{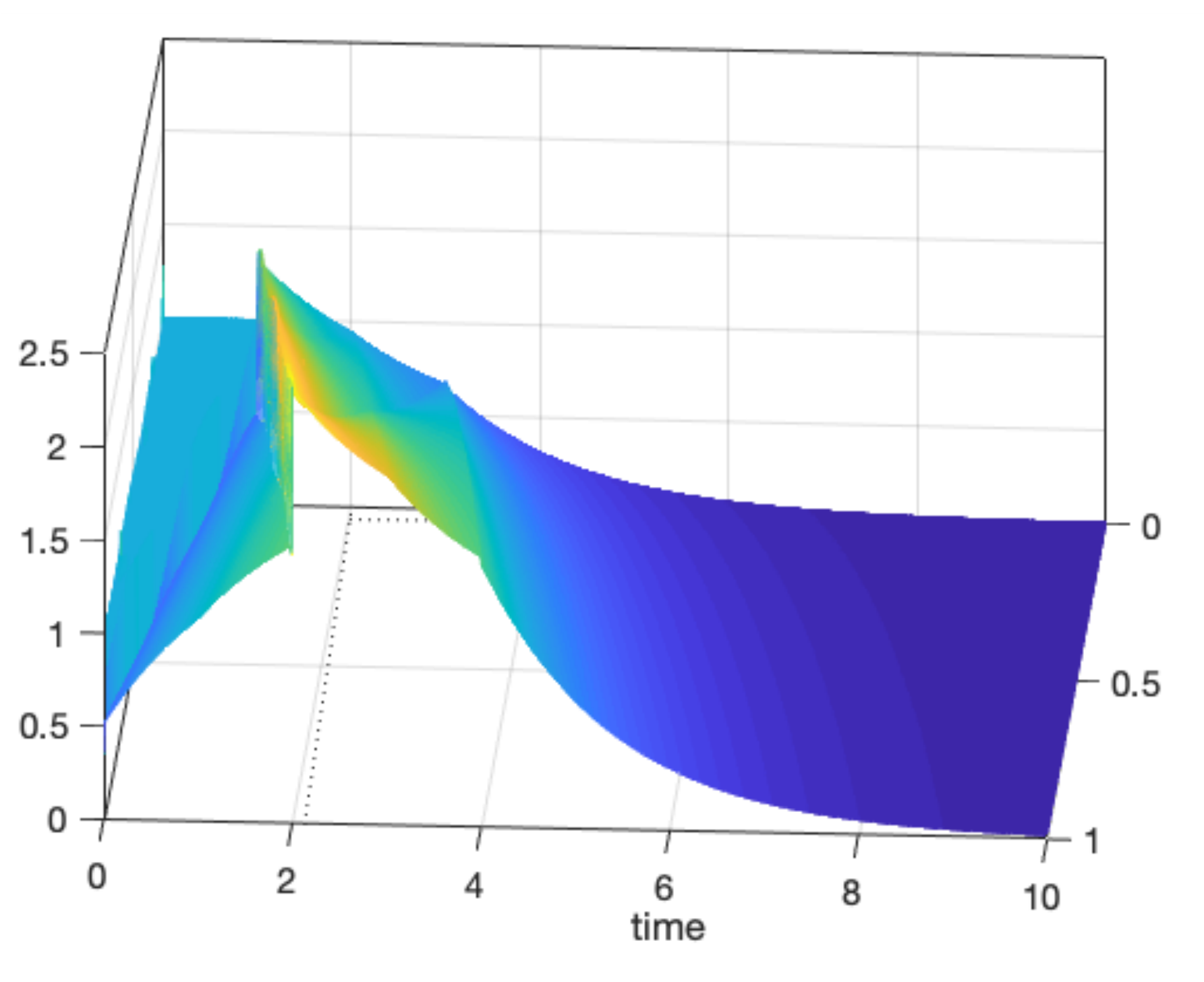}
\caption{Method 1.  $\widetilde{K}$ obtained for nominal $q=3$, but scheduled
$K(q) =K_0+ {\tt feedback}(\widetilde{K},\Phi(q))$. Simulations for $q=2$ top, $q=3$ middle,
$q=4$ bottom \label{method1}}
\end{figure}

\begin{figure}[ht!]
\includegraphics[height=5.2cm]{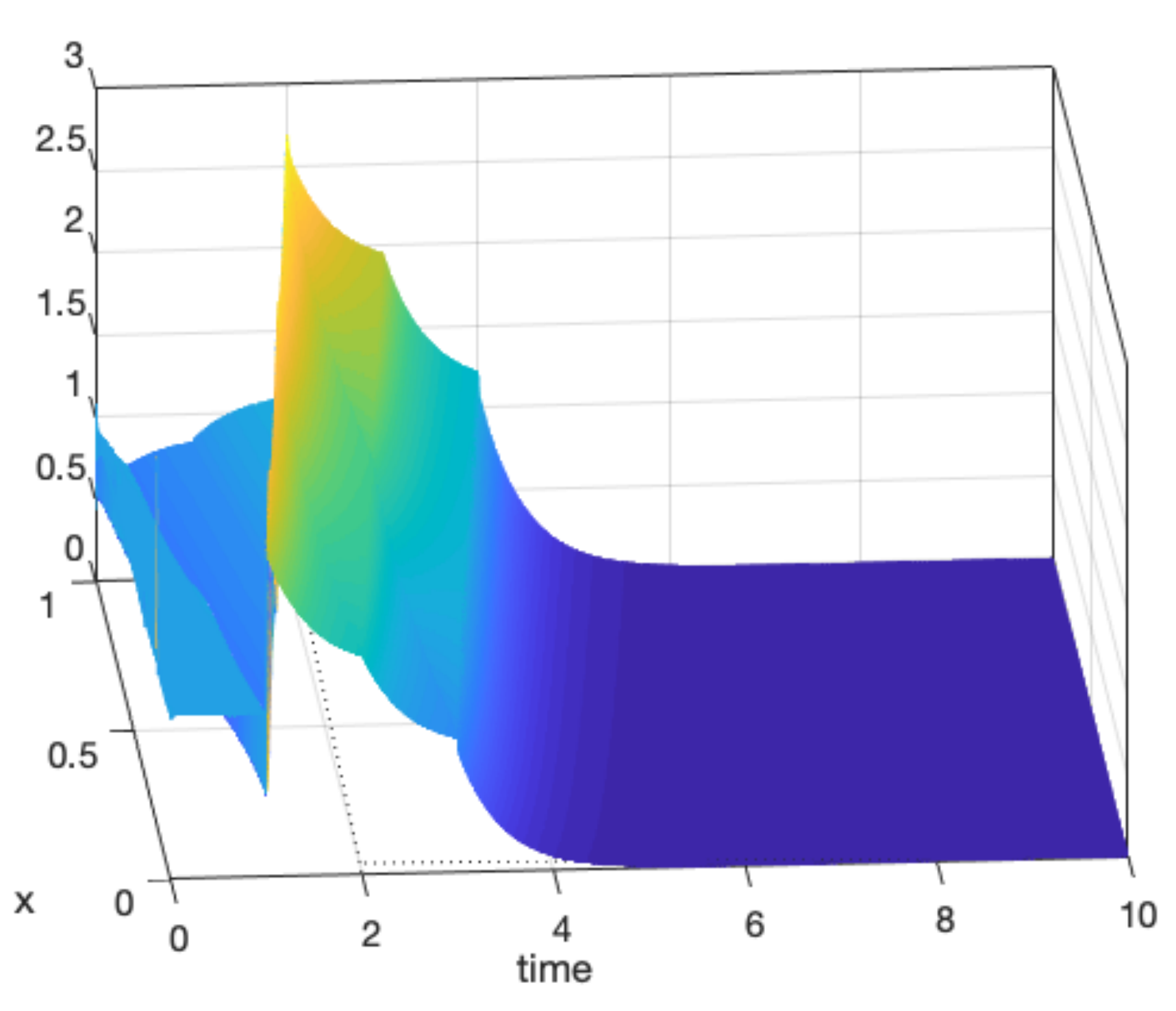}
\includegraphics[height=5.2cm]{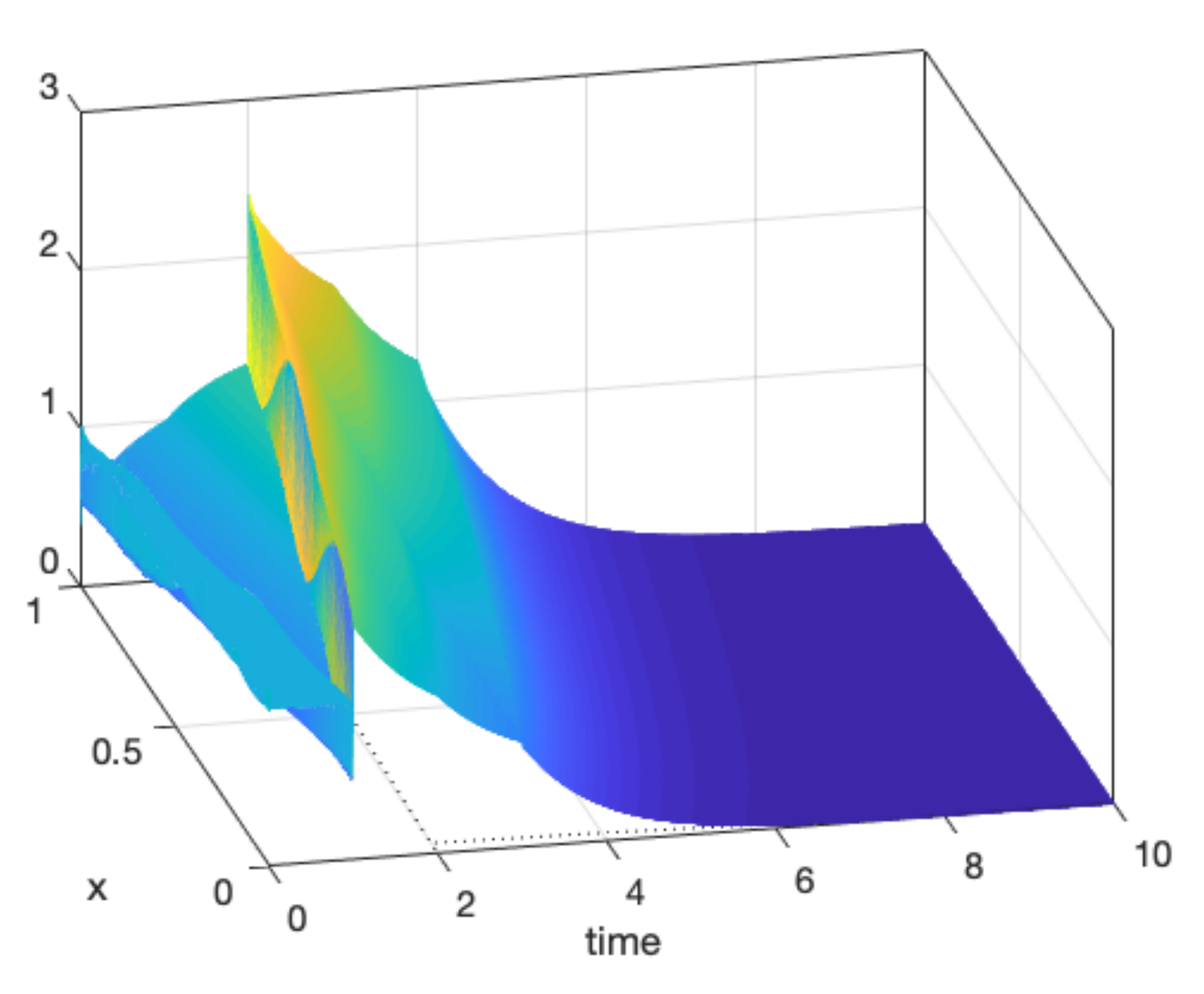}
\includegraphics[height=5.2cm]{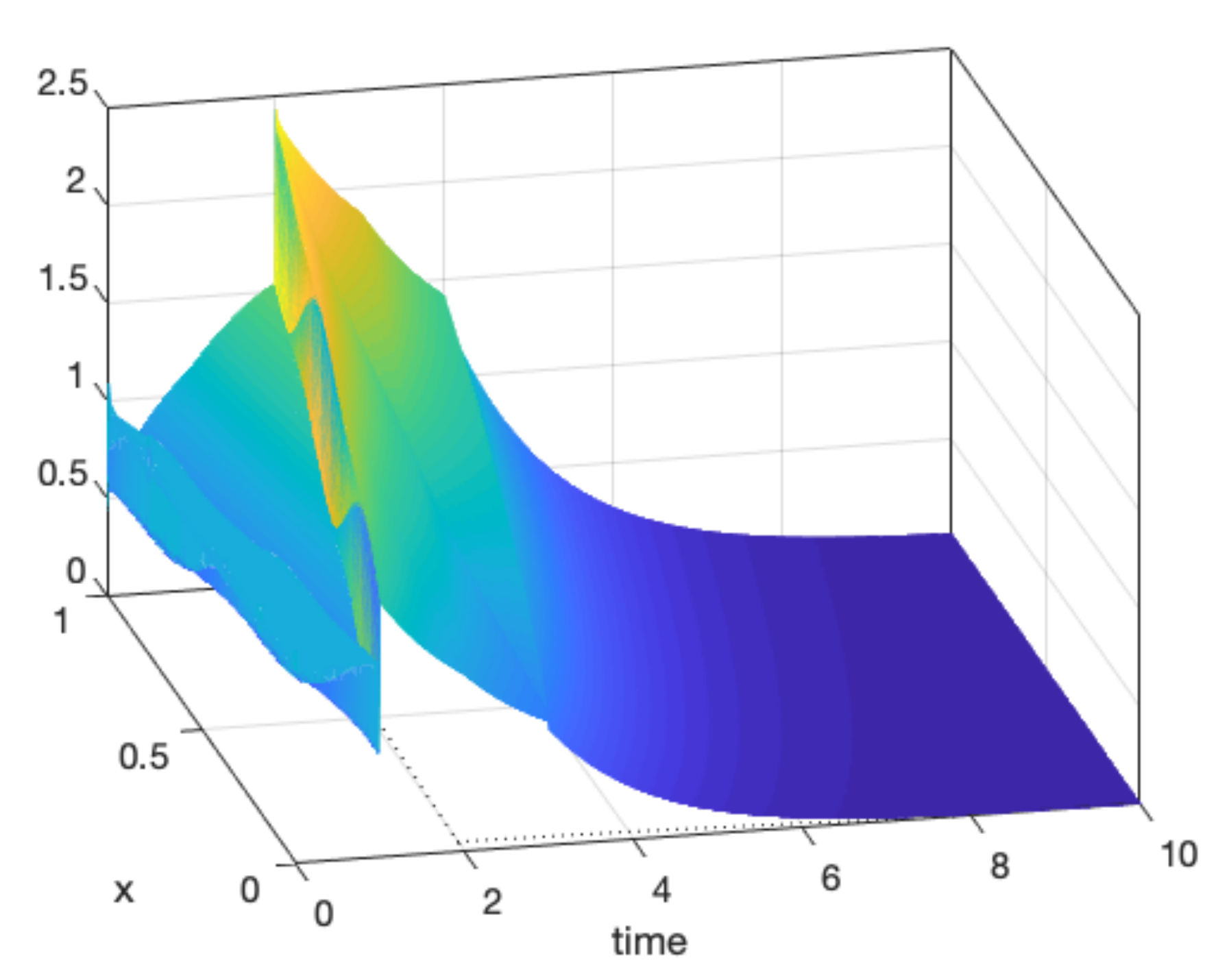}
\caption{Method 2. $\widetilde{K}(q) = \widetilde{K}_{\rm nom} + (q-3)\widetilde{K}_1 + (q-3)^2 \widetilde{K}_2$ and
$K(q) = K_0+{\tt feedback}(\widetilde{K}(q), \Phi(q)$. Simulations for $q=2,3,4$ \label{method2}}
\end{figure}

With the parametric form $\widetilde{K}(q,\x)$ fixed,  we now use again the feedback system $(\widetilde{G}(q),\widetilde{K}(q))$ in
Fig. \ref{swap1} and design a parametric robust controller using the method of 
\cite{AN:2015}, which is included in the {\tt systune} $\,$package and used by default if an
uncertain closed-loop is entered. The tuning goals are chosen as constraints on closed-loop poles
including minimum decay of 0.7, minimum damping of 0.9, with maximum frequency $2$.
The controller obtained is (with $q_0=3$)
\[
\widetilde{K}(q,\x^*) = \widetilde{K}(q_0) + (q-q_0) \widetilde{K}_1(\x^*) + (q-q_0)^2 \widetilde{K}_2(\x^*),
\]
with 
$\widetilde{K}_1 =   [-0.1102, -0.1102, -0.1053]$,  
$\widetilde{K}_2=     [0.03901, 0.03901,  0.02855]$, 
and we retrieve the final parameter varying controller for 
$G(q)$ as
\[
K^{(2)}(q) = K_0 + {\tt feedback}(\widetilde{K}(q,\x^*),-\Phi(q)).
\]
Nominal and scheduled controllers are compared in simulation in Figs. 
\ref{nominal_robust},
\ref{method1},
and \ref{method2}, which
indicate that $K^{(2)}(q)$ achieves the best
performance for frozen-in-time values $q\in [2,4]$.
All controllers are easily implementable, since only real-rational elements in combination with delays are used.

The non-smooth program  (\ref{eq-cast}) was solved with 
{\tt systune}\,   
in $30$s CPU on a Mac OS X 
with $2.66$ GHz Intel Core i7 and $8$ GB RAM. The reader is referred to the MATLAB Control Toolbox 2018b  and higher versions for 
additional examples. More details on this study can be found in \cite{apk_noll:19}.

%




\end{document}